\documentclass[11pt]{article}
\usepackage{latexsym}
\newcommand{\eproof}{\mbox{\ }\hfill $\Box$ \par \vskip 10pt}

\newtheorem{Theorem}{Theorem}[section]
\newtheorem{lemma}[Theorem]{Lemma}
\newtheorem{prop}[Theorem]{Proposition}

\begin{document}

\title{High frequency dispersive estimates for the Schr\"odinger equation in high dimensions}

\author{{\sc Fernando Cardoso, Claudio Cuevas and Georgi Vodev\thanks{Corresponding author}}}

\date{}

\maketitle

\noindent
{\bf Abstract.} We prove optimal dispersive estimates at high frequency for the Schr\"odinger group for a class of real-valued potentials $V(x)=O(\langle x\rangle^{-\delta})$, $\delta>n-1$, and $V\in C^k({\bf R}^n)$, $k>k_n$, where $n\ge 4$ and $\frac{n-3}{2}\le k_n<\frac{n}{2}$.
We also give a sufficient condition in terms of $L^1\to L^\infty$ bounds for the formal iterations of Duhamel's formula, which might be satisfied for potentials of less regularity. 

\setcounter{section}{0}
\section{Introduction and statement of results}

The purpose of this work is to study the question of finding as large as possible class of real-valued potentials $V\in L^\infty({\bf R}^n)$, $n\ge 4$, for which the Schr\"odinger propagator $e^{itG}\chi_a(G)$ satisfies optimal (that is, whitout loss of derivatives) $L^1\to L^\infty$ dispersive estimates, where $G$ denotes the self-adjoint realization of the
operator $-\Delta+V$ on $L^2({\bf R}^n)$, and $\chi_a\in C^\infty({\bf R})$, 
$\chi_a(\lambda)=0$ for $\lambda\le a$, $\chi_a(\lambda)=1$ for $\lambda\ge a+1$, $a\gg 1$. To state our results we need to introduce the class ${\cal C}^k_\delta({\bf R}^n)$, $\delta,k\ge 0$, of all functions $V\in C^k({\bf R}^n)$ satisfying
$$\left\|V\right\|_{{\cal C}_\delta^k}:=\sup_{x\in{\bf R}^n}\sum_{0\le|\alpha|\le k_0}\langle x\rangle^\delta\left|\partial_x^\alpha V(x)\right|$$ $$+\nu\sup_{x\in{\bf R}^n}\sum_{|\beta|=k_0}\langle x\rangle^\delta\sup_{x'\in{\bf R}^n:|x-x'|\le 1}\frac{\left|\partial_x^\beta V(x)-\partial_x^\beta V(x')\right|}{|x-x'|^\nu}<+\infty,$$
where $k_0\ge 0$ is an integer and $\nu=k-k_0$ satisfies $0\le\nu<1$.

\begin{Theorem} Given a $\delta>n-1$, there exists a sequence $\{k_n\}_{n=4}^\infty$, $\frac{n-3}{2}\le k_n<\frac{n}{2}$, so that if $V\in {\cal C}_\delta^k({\bf R}^n)$, $k>k_n$, is a real-valued potential, then we have the following high frequency dispersive estimate
$$\left\|e^{itG}\chi_a(G)\right\|_{L^1\to L^\infty}\le C|t|^{-n/2},\quad t\neq 0,\eqno{(1.1})$$ 
where the constant $C=C(a)>0$ is independent of $t$. 
\end{Theorem}

\noindent
{\bf Remark.} It follows from this theorem and the low frequency dispersive estimates proved in \cite{kn:MV} that if in addition to the assumptions of Theorem 1.1 (or those of Theorem 1.2  below) we assume that zero is neither an eigenvalue nor a resonance of $G$, then we have the following dispersive estimate
$$\left\|e^{itG}P_{ac}\right\|_{L^1\to L^\infty}\le C|t|^{-n/2},\quad t\neq 0,$$ 
where $P_{ac}$ denotes the spectral projection onto the absolutely continuous spectrum of $G$.

Note that with $k_n=\frac{n}{2}$ the above result follows from \cite{kn:MV} where the estimate (1.1) is proved for real-valued potentials $V\in L^\infty({\bf R}^n)$ satisfying 
$$|V(x)|\le C\langle x\rangle^{-\delta},\quad\forall x\in {\bf
R}^n,\eqno{(1.2)}$$ with constants $C>0$, $\delta>n-1$, as well as the condition
$$\widehat V\in L^1.\eqno{(1.3)}$$ 
Previously this has been proved in \cite{kn:JSS} for potentials satisfying (1.2) with $\delta>n$ and (1.3). Proving (1.1) in dimensions $n\ge 4$ without the condition (1.3), however, turns out to be a difficult problem. Note that the potentials in the above theorem do not satisfy (1.3). On the other hand, the counterexample of \cite{kn:GoV} shows that the above theorem cannot hold with $k_n<\frac{n-3}{2}$. Therefore, it is natural to expect that Theorem 1.1 holds with
$k_n=\frac{n-3}{2}$. Indeed, (1.1) has been proved in \cite{kn:CCV1} when $n=4,5$ for potentials $V\in{\cal C}_\delta^k({\bf R}^n)$, $k>\frac{n-3}{2}$, $\delta>3$ if $n=4$, $\delta>5$ if $n=5$. In \cite{kn:CCV2} an analogue of (1.1) with a logarithmic loss of derivatives has been proved for potentials $V\in {\cal C}_\delta^{\frac{n-3}{2}}({\bf R}^n)$, $n=4,5$, where $\delta>3$ if $n=4$, $\delta>5$ if $n=5$. The estimate (1.1) has been recently proved in \cite{kn:EG} when $n=5,7$ for potentials $V$ satisfying (1.2) with $\delta>\frac{3n+5}{2}$ as well as $V\in {\cal C}_\delta^{\frac{n-3}{2}}({\bf R}^n)$, where $\delta>3$ if $n=5$, $\delta>8$ if $n=7$. It also follows from \cite{kn:FY}, \cite{kn:Y} that (1.1) holds for potentials $V$ satisfying (1.2) with $\delta>n+2$ as well as $V\in {\cal C}_\delta^k({\bf R}^n)$ with $\delta>\frac{n}{2}+\frac{2(n-2)}{n-1}$, $k>\frac{n-2}{2}-\frac{1}{n-1}$. 
Note finally that in dimensions one, two and three no regularity of the potential is required in order that (1.1) holds true (see \cite{kn:GS},\cite{kn:S},\cite{kn:M},\cite{kn:RS},\cite{kn:V1},\cite{kn:G}). The same conclusion remains true in dimensions $n\ge 4$ as far as the low and the intermediate frequences are concerned (see \cite{kn:MV}, \cite{kn:V2}).

To prove (1.1) we make use of the semi-classical expansion of the operator $e^{itG}\psi(h^2G)$ obtained in \cite{kn:CV} for potentials satisfying (1.2) with $\delta>\frac{n+2}{2}$, where $\psi\in C_0^\infty((0,+\infty))$ and $0<h\ll 1$. We thus reduce the problem to estimating uniformly in $h$ the $L^1\to L^\infty$ norm of a finite number of operators (denoted by $T_j(t,h)$ below) obtained by iterating the semi-classical Duhamel formula. The advantage is that these operators are defined in terms of the free propagator $e^{itG_0}\psi(h^2G_0)$, where $G_0$ denotes the self-adjoint realization of $-\Delta$ on $L^2({\bf R^n})$ (see Section 2).

In the present paper we also give a sufficient condition for (1.1) to hold in terms of properties of the formal iterations of Duhamel's formula defined as follows (for $t>0$):
$${\cal F}_0(t)=e^{itG_0},\quad {\cal F}_j(t)=i\int_0^t{\cal F}_{j-1}(t-\tau)V{\cal F}_0(\tau)d\tau,\quad j\ge 1.$$
We suppose that there exists a constant
$\varepsilon>0$ such that for all integers $m\ge 1$, $m_1,m_2\ge 0$, we have the bounds
$$\left\|{\cal F}_m(t)\right\|_{L^1\to L^\infty}\le C_mt^{-n/2+\varepsilon m},\quad 0<t\le 1,\eqno{(1.4)}$$
$$\left\|\int_{I(\gamma)}{\cal F}_{m_1}(t-\tau)V{\cal F}_{m_2}(\tau)d\tau\right\|_{L^1\to L^\infty}\le C_{m_1,m_2}\gamma^{\varepsilon}t^{-n/2},\quad\forall t>0,\eqno{(1.5)}$$
where $0<\gamma\le 1$, $I(\gamma)\subset[0,t]$ is an interval either of the form $[0,\gamma_1]$ or of the form $[t-\gamma_1,t]$, $\gamma_1=t/2$ if $t\le 2\gamma$, $\gamma_1=\gamma$ if $t\ge 2\gamma$. 

\begin{Theorem} Let $V$ satisfy (1.2) with $\delta>n$ and suppose (1.4) and (1.5) fulfilled. Then, the dispersive estimate (1.1) holds true for all $t>0$.
\end{Theorem}

It is easy to see that if $V$ satisfies (1.3), then (1.4) and (1.5) hold with $\varepsilon=1$. However, it might happend that (1.4) and (1.5) hold true for potentials of less regularity. In fact, we expect that (1.4) and (1.5) hold for potentials $V\in {\cal C}_\delta^k({\bf R}^n)$ with $\delta>\frac{n+1}{2}$, $k>\frac{n-3}{2}$. Indeed, this has been proved in \cite{kn:CCV1} for $m=1$, $m_1=m_2=0$. The problem, however, gets much harder for $m\ge 2$, $m_1,m_2\ge 1$.

To prove Theorem 1.2 we take advantage of the analysis carried out in \cite{kn:V2} under the only assumption that $V$ satisfies (1.2) with $\delta>\frac{n+2}{2}$ (see Section 5).

\section{Reduction to semi-classical dispersive estimates}

Set 
$$F(t)=i\int_0^t e^{i(t-\tau)G_0}Ve^{i\tau G_0}d\tau,\quad t>0.$$
It is easy to see that (1.1) is a consequence of the following

\begin{Theorem} Under the assumptions of Theorem 1.1, the following dispersive estimates hold true for all $0<h\ll 1$, $t>0$:
$$\|F(t)\|_{L^1\to L^\infty}\le Ct^{-n/2},\eqno{(2.1)}$$
$$\left\|e^{itG}\psi(h^2G)-e^{itG_0}\psi(h^2G_0)-F(t)\psi(h^2G_0)\right\|_{L^1\to L^\infty}\le Ch^\beta t^{-n/2},\eqno{(2.2)}$$
with some constants $C,\beta>0$ independent of $t$ and $h$.
\end{Theorem}

The estimate (2.1) is proved in \cite{kn:CCV1} for potentials $V\in {\cal C}_\delta^k({\bf R}^n)$, $\delta>n-1$, $k>\frac{n-3}{2}$. In what follows we will derive (2.2) from the semi-classical expansion obtained in \cite{kn:CV} and based on the following semi-classical version of Duhamel's formula
$$e^{itG}\psi(h^2G)=Q(h)e^{itG_0}\psi_1(h^2G_0)\psi(h^2G)+i\int_0^tQ(h)\psi_1(h^2G_0)
e^{i(t-\tau)G_0}Ve^{i\tau G}\psi(h^2G)d\tau,\eqno{(2.3)}$$
where $\psi_1\in C_0^\infty((0,+\infty))$, $\psi_1=1$ on supp$\,\psi$, and 
$$Q(h)=\left(1+\psi_1(h^2G_0)-\psi_1(h^2G)\right)^{-1}.$$
Iterating (2.3) $m$ times we get
$$e^{itG}\psi(h^2G)=\sum_{j=0}^m T_j(t,h)+\int_0^t R_{m}(t-\tau,h)e^{i\tau G}\psi(h^2G)d\tau,\eqno{(2.4)}$$
$$T_0(t,h)=Q(h)e^{itG_0}\psi_1(h^2G_0)\psi(h^2G),$$
$$T_j(t,h)=\int_0^t R_{j-1}(t-\tau,h)T_0(\tau,h)d\tau,\quad j\ge 1,$$
where the operators $R_j$ are defined as follows
$$R_0(t,h)=iQ(h)e^{itG_0}\psi_1(h^2G_0)V,$$
$$R_j(t,h)=\int_0^t R_{j-1}(t-\tau,h)R_0(\tau,h)d\tau,\quad j\ge 1.$$
The following dispersive estimates are proved in \cite{kn:CV} (see Theorem 1.3) (it is easy to see that the $\epsilon$ there can be taken zero).

\begin{prop} Assume that $V$ satisfies (1.2) with $\delta>\frac{n+2}{2}$. Then the following dispersive estimates hold true for all $t>0$, $0<h\ll 1$,
$$\left\|T_j(t,h)\right\|_{L^1\to L^\infty}\le C_jh^{j-n/2}t^{-n/2},\quad j\ge 1,\eqno{(2.5)}$$
$$\left\|e^{itG}\psi(h^2G)-\sum_{j=0}^m T_j(t,h)\right\|_{L^1\to L^\infty}\le C_mh^{m+1-n/2} t^{-n/2},\quad m\ge 1.\eqno{(2.6)}$$
\end{prop}

We have (e.g. see Lemma A.1 of \cite{kn:MV})
$$\psi_1(h^2G_0)-\psi_1(h^2G)=O(h^2): L^1\to L^1,$$
so
$$Q(h)=Id+O(h^2): L^1\to L^1.$$
Therefore,
$$\left\|T_0(t,h)-e^{itG_0}\psi(h^2G_0)\right\|_{L^1\to L^\infty}\le Ch^{2} t^{-n/2}.\eqno{(2.7)}$$
Clearly, the estimate (2.2) follows from combining Proposition 2.2, (2.7) and the following

\begin{prop} Given a $\delta>n-1$, there exists a sequence $\{k_n\}_{n=4}^\infty$, $\frac{n-3}{2}\le k_n<\frac{n}{2}$, so that if $V\in{\cal C}_\delta^k({\bf R}^n)$, $k>k_n$, then we have the estimates
$$\left\|T_1(t,h)-F(t)\psi(h^2G_0)\right\|_{L^1\to L^\infty}\le Ch^\beta t^{-n/2},\eqno{(2.8)}$$
$$\left\|T_j(t,h)\right\|_{L^1\to L^\infty}\le Ch^\beta t^{-n/2},\quad 2\le j\le n/2,\eqno{(2.9})$$
with some constants $C,\beta>0$ independent of $h$ and $t$.
\end{prop}

\section{Proof of Proposition 2.3}

Set
$$\widetilde T_1(t,h)=T_1(t,h)-i\int_0^t e^{i(t-\tau)G_0}\psi_1(h^2G_0)Ve^{i\tau G_0}\psi(h^2G_0)d\tau.$$
It is proved in \cite{kn:CCV1} (see Proposition 2.6) that if $V\in{\cal C}_\delta^k({\bf R}^n)$
with $\delta>n-1$, $k>(n-3)/2$, then
$$\left\|i\int_0^t e^{i(t-\tau)G_0}\psi_1(h^2G_0)Ve^{i\tau G_0}\psi(h^2G_0)d\tau-F(t)\psi(h^2G_0)\right\|_{L^1\to L^\infty}\le Ch^\beta t^{-n/2}\eqno{(3.1)}$$
with constants $C,\beta>0$ independent of $h$ and $t$. Therefore, to prove (2.8) it suffices to show that the operator $\widetilde T_1$ satisfies the estimate
$$\left\|\widetilde T_1(t,h)\right\|_{L^1\to L^\infty}\le Ch^\beta t^{-n/2}\eqno{(3.2)}$$
with constants $C,\beta>0$ independent of $h$ and $t$. It is easy also to see that the operators $T_j$, $j\ge 2$, are of the form $Q(h)\widetilde T_j(t,h)\psi(h^2G)$. Therefore, it suffices to prove (2.9) with $T_j$ replaced by $\widetilde T_j$.

Let $\rho\in C_0^\infty({\bf R}^n)$, $\rho\ge 0$, be a real-valued function such that $\int \rho(x)dx=1$, and set $\rho_\theta(x)=\theta^{-n}\rho(x/\theta)$, where $0<\theta\le 1$. Let $V\in{\cal C}_\delta^k({\bf R}^n)$ with $\delta>n-1$, where $k$ will be fixed later on such that 
$\frac{n-3}{2}< k<\frac{n}{2}$. Set $V_\theta=V*\rho_\theta$. It is easy to see that we have the bounds
$$\left|V_\theta(x)\right|\le C\langle x\rangle^{-\delta},\quad\forall x\in {\bf
R}^n,\eqno{(3.3)}$$
$$\left|V_\theta(x)-V(x)\right|\le C\theta^{1/2}\langle x\rangle^{-\delta},\quad\forall x\in {\bf R}^n,\eqno{(3.4)}$$
$$\left|\partial_x^\alpha V_\theta(x)\right|\le C\langle x\rangle^{-\delta},\quad\forall x\in {\bf R}^n,\,|\alpha|\le k_0,\eqno{(3.5)}$$
$$\left|\partial_x^\alpha V_\theta(x)\right|\le C_\alpha \theta^{k-|\alpha|}\langle x\rangle^{-\delta},\quad\forall x\in {\bf R}^n,\,|\alpha|\ge k_0+1,\eqno{(3.6)}$$
where $k-1<k_0\le k$ is an integer. Let us also see that 
$$\left\|\widehat V_\theta\right\|_{L^1}\le C_\epsilon\theta^{k-n/2-\epsilon},\quad\forall\, 0<\epsilon\ll 1.\eqno{(3.7)}$$
Since ${\cal C}_\delta^k({\bf R}^n)\subset H^{k-\epsilon/2}({\bf R}^n)$, $\forall\,0<\epsilon\ll 1$, we have $\langle\xi\rangle^{k-\epsilon/2}\widehat V(\xi)\in L^2({\bf R}^n)$. Hence $\widehat V\in L^p({\bf R}^n)$,
where $\frac{1}{p}=\frac{1}{2}+\frac{k}{n+\epsilon}$. We have 
$$\left\|\widehat V_\theta\right\|_{L^1}=\left\|\widehat V\widehat\rho_\theta\right\|_{L^1}\le
\left\|\widehat V\right\|_{L^p}\left\|\widehat\rho_\theta\right\|_{L^q}=C\theta^{-n/q},$$
where $\frac{1}{q}=\frac{1}{2}-\frac{k}{n+\epsilon}$, which clearly implies (3.7).

Let $G_\theta$ denote the self-adjoint realization of $-\Delta+V_\theta(x)$ on $L^2({\bf R}^n)$.
Denote also by $Q_\theta(h)$ the operator obtained by replacing in the definition of $Q(h)$ the operator $G$ by $G_\theta$. Define the operators $\widetilde T_{j,\theta}$  by replacing $Q(h)$ and $V$ by $Q_\theta(h)$ and $V_\theta$, respectively, in the definition of $\widetilde T_{j}$. In the case of $\widetilde T_1$ we replace only those $V$ and $Q(h)$ staying between the operators $e^{i(t-\tau)G_0}$ and $e^{i\tau G_0}$. Using (3.3) and (3.4) we will prove the following

\begin{prop} The following dispersive estimates hold for all $t>0$, $0<h\ll 1$, $0<\theta\le 1$,
$$\left\|\widetilde T_j(t,h)-\widetilde T_{j,\theta}(t,h)\right\|_{L^1\to L^\infty}\le C\theta^{1/2}h^{j-n/2} t^{-n/2},\quad 1\le j\le n/2.\eqno{(3.8)}$$
\end{prop}

{\it Proof.} We write
$$\widetilde T_1(t,h)-\widetilde T_{1,\theta}(t,h)$$ $$=iQ(h)\int_0^te^{i(t-\tau)G_0}\psi_1(h^2G_0)\left(VQ(h)-V_\theta Q_\theta(h)\right)e^{i\tau G_0}\psi_1(h^2G_0)d\tau\left(\psi(h^2G)-\psi(h^2G_0)\right)$$
$$+i\left(Q(h)-1\right)\int_0^te^{i(t-\tau)G_0}\psi_1(h^2G_0)\left(VQ(h)-V_\theta Q_\theta(h)\right)e^{i\tau G_0}\psi(h^2G_0)d\tau$$
$$+i\int_0^te^{i(t-\tau)G_0}\psi_1(h^2G_0)\left(V(Q(h)-1)-V_\theta( Q_\theta(h)-1)\right)e^{i\tau G_0}\psi(h^2G_0)d\tau$$
$$=:\sum_{j=1}^3{\cal P}^{(j)}_{\theta}(t,h).\eqno{(3.9)}$$
Define the operators $F_j(t,h)$, $F_{j,\theta}(t,h)$, $j=0,1,...,$ by
$$F_0(t,h)=F_{0,\theta}(t,h)=e^{itG_0}\psi_1(h^2G_0),$$
$$F_j(t,h)=i\int_0^tF_0(t-\tau,h)VQ(h)F_{j-1}(\tau,h)d\tau,\quad j\ge 1,$$
$$F_{j,\theta}(t,h)=i\int_0^tF_0(t-\tau,h)V_\theta Q_\theta(h)F_{j-1,\theta}(\tau,h)d\tau,\quad j\ge 1.$$
Clearly, $\widetilde T_j=F_j$, $\widetilde T_{j,\theta}=F_{j,\theta}$ for $j\ge 2$. We write
$$F_j(t,h)-F_{j,\theta}(t,h)=i\int_0^t F_0(t-\tau,h)\left(VQ(h)-V_\theta Q_\theta(h)\right) F_{j-1,\theta}(\tau,h)d\tau$$ $$+i\int_0^t F_0(t-\tau,h)VQ(h)\left( F_{j-1}(\tau,h)-F_{j-1,\theta}(\tau,h)\right)d\tau.\eqno{(3.10)}$$
Let us see that (3.8) follows from the following estimates

\begin{prop} For all $t>0$, $0<h\ll 1$, $0<\theta\le 1$, $1/2-\epsilon/2\le s\le (n-1)/2$, $0<\epsilon\ll 1$, $j\ge 0$, we have the estimates
$$\left\|\langle x\rangle^{-1/2-s-\epsilon}F_j(t,h)\right\|_{L^1\to L^2}\le C_jh^{j+s-(n-1)/2}t^{-s-1/2},\eqno{(3.11)}$$
$$\left\|\langle x\rangle^{-1/2-s-\epsilon}\left(F_j(t,h)-F_{j,\theta}(t,h)\right)\right\|_{L^1\to L^2}\le C_j\theta^{1/2}h^{j+s-(n-1)/2}t^{-s-1/2}.\eqno{(3.12)}$$
\end{prop}

\noindent
{\bf Remark.} The estimate (3.11) with $j=0$ holds true for all $t\neq 0$. In other words, the adjoint of the operator 
$${\cal A}=F_0(t,h)\langle x\rangle^{-1/2-s-\epsilon}:L^2\to L^\infty$$
satisfies (3.11) with $j=0$, and hence so does ${\cal A}$. This will be often used below.

We need the following

\begin{lemma} For all $0<h\le h_0$, $0<\theta\le 1$, $0\le s\le\delta$, we have the bounds
$$\left\|\langle x\rangle^{-s}Q(h)\langle x\rangle^{s}\right\|_{L^2\to L^2}\le C,\eqno{(3.13)}$$
$$\left\|\langle x\rangle^{-s}\left(Q(h)-Q_\theta(h)\right)\langle x\rangle^{s}\right\|_{L^2\to L^2}\le C\theta^{1/2},\eqno{(3.14)}$$
with constants $C,h_0>0$ independent of $h$ and $\theta$.
\end{lemma}

{\it Proof.} Clearly, (3.13) follows from the bound
$$\left\|\langle x\rangle^{-s}\left(\psi_1(h^2G)-\psi_1(h^2G_0)\right)\langle x\rangle^{s}\right\|_{L^2\to L^2}\le Ch^2,\eqno{(3.15)}$$
 proved in \cite{kn:V3} (see Lemma 2.3). To prove (3.14) we write
 $$Q(h)-Q_\theta(h)=\left(\psi_1(h^2G_\theta)-\psi_1(h^2G)\right)Q(h)+\left(\psi_1(h^2G_\theta)-\psi_1(h^2G_0)\right)(Q(h)-Q_\theta(h)).$$
Therefore, (3.14) follows from combining (3.13), (3.15) and the bound
$$\left\|\langle x\rangle^{-s}\left(\psi_1(h^2G)-\psi_1(h^2G_\theta)\right)\langle x\rangle^{s}\right\|_{L^2\to L^2}\le C\theta^{1/2}h^2.\eqno{(3.16)}$$
To prove (3.16) we will use the Helffer-Sj\"ostrand formula
$$\psi_1(h^2G)=\frac{2}{\pi}\int_{\bf C}\frac{\partial\widetilde\varphi}{\partial\bar z}(z)(h^2G-z^2)^{-1}zL(dz),\eqno{(3.17)}$$
where $L(dz)$ denotes the Lebesgue measure on ${\bf C}$, $\widetilde\varphi\in C_0^\infty({\bf C})$ is an almost analytic continuation of
$\varphi(\lambda)=\psi_1(\lambda^2)$, supported in a small complex neighbourhood of supp$\,\varphi$ and satisfying
$$\left|\frac{\partial\widetilde\varphi}{\partial\bar z}(z)\right|\le C_N|{\rm Im}\,z|^N,\quad\forall N\ge 1.$$
In view of (3.17) we can write
$$\psi_1(h^2G)-\psi_1(h^2G_\theta)=\frac{2h^2}{\pi}\int_{\bf C}\frac{\partial\widetilde\varphi}{\partial\bar z}(z)(h^2G_\theta-z^2)^{-1}(V_\theta-V)(h^2G-z^2)^{-1}zL(dz). \eqno{(3.18)}$$
It is shown in \cite{kn:V3} (see the proof of Lemma 2.3) that the free resolvent satisfies the bound (for $z\in{\rm supp}\,\widetilde\varphi$)
 $$\left\|\langle x\rangle^{-s}(h^2G_0-z^2)^{-1}\langle x\rangle^{s}\right\|_{L^2\to L^2}\le C_1|{\rm Im}\,z|^{-q},\quad{\rm Im}\,z\neq 0,\eqno{(3.19)}$$
with constants $C_1,q>0$ independent of $z$ and $h$. By (3.19) and the identity
$$(h^2G-z^2)^{-1}=(h^2G_0-z^2)^{-1}-h^2(h^2G-z^2)^{-1}V(h^2G_0-z^2)^{-1},$$
we obtain (for $z\in{\rm supp}\,\widetilde\varphi$, $0\le s\le\delta$)
$$\left\|\langle x\rangle^{-s}(h^2G-z^2)^{-1}\langle x\rangle^{s}\right\|_{L^2\to L^2}\le \left\|\langle x\rangle^{-s}(h^2G_0-z^2)^{-1}\langle x\rangle^{s}\right\|_{L^2\to L^2}$$ $$+Ch^2\left\|(h^2G-z^2)^{-1}\right\|_{L^2\to L^2}\left\|\langle x\rangle^{-s}(h^2G_0-z^2)^{-1}\langle x\rangle^{s}\right\|_{L^2\to L^2}$$
$$\le C_2|{\rm Im}\,z|^{-q-1},\quad{\rm Im}\,z\neq 0.\eqno{(3.20)}$$
By (3.4), (3.18) and (3.20),
$$\left\|\langle x\rangle^{-s}\left(\psi_1(h^2G)-\psi_1(h^2G_\theta)\right)\langle x\rangle^{s}\right\|_{L^2\to L^2}$$ $$\le C\theta^{1/2}h^2\int_{\bf C}\left|\frac{\partial\widetilde\varphi}{\partial\bar z}(z)\right|\left\|(h^2G_\theta-z^2)^{-1}\right\|_{L^2\to L^2}\left\|\langle x\rangle^{-s}(h^2G-z^2)^{-1}\langle x\rangle^{s}\right\|_{L^2\to L^2}L(dz)$$
$$\le C\theta^{1/2}h^2\int_{\bf C}\left|\frac{\partial\widetilde\varphi}{\partial\bar z}(z)\right||{\rm Im}\,z|^{-q-2}L(dz)\le C\theta^{1/2}h^2.\eqno{(3.21)}$$
\eproof

Using (3.3), (3.4), Lemma 3.3 and (3.11) with $j=0$, we obtain
$$\left\|{\cal P}^{(3)}_\theta(t,h)\right\|_{L^1\to L^\infty}$$ $$\le C\theta^{1/2}\int_0^{t/2}\left\|e^{i(t-\tau)G_0}\psi_1(h^2G_0)\langle x\rangle^{-n/2-\epsilon}\right\|_{L^2\to L^\infty}\left\|\langle x\rangle^{-1-\epsilon}e^{i\tau G_0}\psi(h^2G_0)\right\|_{L^1\to L^2}d\tau$$
$$+C\theta^{1/2}\int_{t/2}^t\left\|e^{i(t-\tau)G_0}\psi_1(h^2G_0)\langle x\rangle^{-1-\epsilon}\right\|_{L^2\to L^\infty}\left\|\langle x\rangle^{-n/2-\epsilon}e^{i\tau G_0}\psi(h^2G_0)\right\|_{L^1\to L^2}d\tau$$
 $$\le C\theta^{1/2}t^{-n/2}h^{1-n/2-\epsilon/2}\int_0^h\tau^{-1+\epsilon/2}d\tau+
 C\theta^{1/2}t^{-n/2}h^{1-n/2+\epsilon/2}\int_h^\infty\tau^{-1-\epsilon/2}d\tau$$ $$\le C\theta^{1/2}t^{-n/2}h^{1-n/2}.$$
 Clearly, the $L^1\to L^\infty$ norm of the operators ${\cal P}^{(j)}_\theta$, $j=1,2$, can be bounded in the same way. Let now $j\ge 2$. Using (3.10), Proposition 3.2 and Lemma 3.3, we obtain
 
 $$\left\|F_j(t,h)-F_{j,\theta}(t,h)\right\|_{L^1\to L^\infty}$$ $$\le C\theta^{1/2}\int_0^{t/2}\left\|F_0(t-\tau,h)\langle x\rangle^{-n/2-\epsilon}\right\|_{L^2\to L^\infty}\left\|\langle x\rangle^{-1-\epsilon}F_{j-1,\theta}(t,h)\right\|_{L^1\to L^2}d\tau$$
$$+C\theta^{1/2}\int_{t/2}^t\left\|F_0(t-\tau,h)\langle x\rangle^{-1-\epsilon}\right\|_{L^2\to L^\infty}\left\|\langle x\rangle^{-n/2-\epsilon}F_{j-1,\theta}(t,h)\right\|_{L^1\to L^2}d\tau$$
 $$+ C\int_0^{t/2}\left\|F_0(t-\tau,h)\langle x\rangle^{-n/2-\epsilon}\right\|_{L^2\to L^\infty}\left\|\langle x\rangle^{-1-\epsilon}\left(F_{j-1}(\tau,h)-F_{j-1,\theta}(\tau,h)\right)\right\|_{L^1\to L^2}d\tau$$
$$+C\int_{t/2}^t\left\|F_0(t-\tau,h)\langle x\rangle^{-1-\epsilon}\right\|_{L^2\to L^\infty}\left\|\langle x\rangle^{-n/2-\epsilon}\left(F_{j-1}(\tau,h)-F_{j-1,\theta}(\tau,h)\right)\right\|_{L^1\to L^2}d\tau$$
 $$\le C\theta^{1/2}t^{-n/2}h^{j-n/2-\epsilon/2}\int_0^h\tau^{-1+\epsilon/2}d\tau+
 C\theta^{1/2}t^{-n/2}h^{j-n/2+\epsilon/2}\int_h^\infty\tau^{-1-\epsilon/2}d\tau$$ $$\le C\theta^{1/2}t^{-n/2}h^{j-n/2}.$$
 \eproof
 
 {\it Proof of Proposition 3.2.} The estimate (3.11) with $j=0$ is proved in \cite{kn:V2} (see (2.1)). By induction in $j$, it is easy to see that (3.11) for any $j$ follows from this and the following well-known estimate
 $$\left\|\langle x\rangle^{-s}e^{itG_0}\psi_1(h^2G_0)\langle x\rangle^{-s}\right\|_{L^2\to L^2}\le C\langle t/h\rangle^{-s},\quad s\ge 0.\eqno{(3.22)}$$
 Similarly, using (3.10) together with Lemma 3.3, one can easily get (3.12).
 \eproof
 
 \section{Study of the operators $\widetilde T_{j,\theta}$}
 
 We will first show that the estimates (2.9) and (3.2) follow from Proposition 3.1 and the following
 
 \begin{prop} Let $V\in {\cal C}_\delta^k({\bf R}^n)$ with $\delta>n-1$, $\frac{n-3}{2}<k<\frac{n}{2}$. Then, there exist a constant $\varepsilon_0>0$ and a sequence $\{p_j\}_{j=1}^\infty$, $p_j>0$, depending on $\delta$ but independent of $k$, so that for all $0<h\ll 1$, $0<\theta\le 1$, 
 $0<\epsilon\ll 1$, $t>0$, $j\ge 1$, satisfying
  $$h^2\theta^{k-n/2-\epsilon}\ll 1,\eqno{(4.1)}$$
  we have the estimate
 $$\left\|\widetilde T_{j,\theta}(t,h)\right\|_{L^1\to L^\infty}\le C_jh^{\varepsilon_0}t^{-n/2}+
 C_{j,\epsilon}h^{p_j}\theta^{-j(n/2-k+\epsilon)}t^{-n/2},\eqno{(4.2)}$$
 where $C_j,C_{j,\epsilon}>0$ are independent of $t$, $h$ and $\theta$.
 \end{prop}
 
 Fix an integer $1\le j\le n/2$. Take $\theta=h^{n+1-2j}$ and set
 $$k_n^{(j)}=\frac{n}{2}-\frac{\min\left\{\frac{3}{2},\frac{p_j}{j}\right\}}{n+1-2j}.$$
 It is easy to see that if $k_n^{(j)}<k<n/2$ and $\epsilon$ is taken small enough, we can arrange (4.1), and the estimates (3.8) and (4.2) imply
  $$\left\|\widetilde T_j(t,h)\right\|_{L^1\to L^\infty}\le Ch^\beta t^{-n/2},\eqno{(4.3)}$$
  with some $C,\beta>0$. Thus, taking
  $$k_n=\max_{1\le j\le n/2}k_n^{(j)}$$
  we get the desired result.\\

 {\it Proof of Proposition 4.1.} We need the following
 
 \begin{lemma} For all $0<\theta\le 1$, $0<\epsilon\ll 1$, $t\in{\bf R}$, we have the estimate
 $$\left\|e^{-itG_0}V_\theta e^{itG_0}\right\|_{L^1\to L^1}\le C_\epsilon \theta^{k-n/2-\epsilon},\eqno{(4.4)}$$
 with a constant $C_\epsilon>0$ independent of $t$ and $\theta$. Moreover, given any integer $m\ge 1$, the operator $Q_\theta(h)$ can be decomposed as $P_m^{(1)}(h,\theta)+P_m^{(2)}(h,\theta)$, where the operator $P_m^{(1)}$ satisfies the estimate
 $$\left\|e^{-itG_0}P_m^{(1)}(h,\theta) e^{itG_0}\right\|_{L^1\to L^1}\le 2,\eqno{(4.5)}$$
 for all $t\in{\bf R}$ and all $0<h\ll 1$, $0<\theta\le 1$ such that (4.1) holds, while
 the operator $P_m^{(2)}$ satisfies the estimate
 $$\left\|\langle x\rangle^{-s}P_m^{(2)}(h,\theta)\langle x\rangle^{s}\right\|_{L^2\to L^2}\le C_mh^{2m+2},\quad 0\le s\le\delta,\eqno{(4.6)}$$
where $C_m>0$ is independent of $h$ and $\theta$.
 \end{lemma}
 
 {\it Proof.} The estimate (4.4) follows from (3.7) and the following estimate proved in \cite{kn:JSS}:
 $$\left\|e^{-itG_0}V_\theta e^{itG_0}\right\|_{L^1\to L^1}\le \left\|\widehat V_\theta\right\|_{L^1}.\eqno{(4.7)}$$
 To decompose the operator $Q_\theta(h)$ we will use the formula (3.17) together with the resolvent identity
 $$\left(h^2G_\theta-z^2\right)^{-1}-\left(h^2G_0 -z^2\right)^{-1} =\sum_{j=1}^m\left(h^2G_0-z^2\right)^{-1}\left(-h^2V_\theta\left(h^2G_0-z^2\right)^{-1}\right)^j$$ $$+\left(h^2G_\theta-z^2\right)^{-1}\left(-h^2V_\theta\left(h^2G_0-z^2\right)^{-1}\right)^{m+1}:=\sum_{\ell=1}^2{\cal M}_m^{(\ell)}(z,h,\theta).$$
 Set
 $$M_m^{(\ell)}(h,\theta)=\frac{2}{\pi}\int_{\bf C}\frac{\partial\widetilde\varphi}{\partial\bar z}(z){\cal M}_m^{(\ell)}(z,h,\theta)zL(dz),$$
  $$P_m^{(1)}(h,\theta)=\left(1-M_m^{(1)}(h,\theta)\right)^{-1},$$
  $$P_m^{(2)}(h,\theta)=\left(1-M_m^{(1)}(h,\theta)-M_m^{(2)}(h,\theta)\right)^{-1}-\left(1-M_m^{(1)}(h,\theta)\right)^{-1}$$
  $$=\left(1-M_m^{(1)}(h,\theta)-M_m^{(2)}(h,\theta)\right)^{-1}M_m^{(2)}(h,\theta)\left(1-M_m^{(1)}(h,\theta)\right)^{-1}.$$
 We need now the following well known estimate (e.g. see (2.14) of \cite{kn:MV}):
 $$\left\|\left(h^2G_0 -z^2\right)^{-1} \right\|_{L^1\to L^1}\le C\left|{\rm Im}\,z\right|^{-q},\eqno{(4.8)}$$
 for $z\in{\rm supp}\,\widetilde\varphi$, ${\rm Im}\,z\neq 0$, $0<h\le 1$,
 where the constants $C,q>0$ are independent of $z$ and $h$. By (4.7) and (4.8), we get
 $$\left\|e^{-itG_0}{\cal M}_m^{(1)}(z,h,\theta) e^{itG_0}\right\|_{L^1\to L^1}\le\sum_{j=1}^m  C_j\left\|h^2\widehat V_\theta\right\|_{L^1}^j|{\rm Im}\,z|^{-q(j+1)}.\eqno{(4.9)}$$
 By (3.7) and (4.9), we conclude
 $$\left\|e^{-itG_0}M_m^{(1)}(h,\theta) e^{itG_0}\right\|_{L^1\to L^1}\le C_m\left\|h^2\widehat V_\theta\right\|_{L^1}\left(1+\left\|h^2\widehat V_\theta\right\|_{L^1}\right)^{m-1}$$ $$\le\widetilde C_mh^2\theta^{k-n/2-\epsilon}\le  1/2,\eqno{(4.10)}$$
 provided (4.1) is satisfied. Clearly, (4.5) follows from (4.10). On the other hand, it is easy to see that (4.6) follows from the estimates
 $$\left\|\langle x\rangle^{-s}M_m^{(1)}(h,\theta)\langle x\rangle^{s}\right\|_{L^2\to L^2}\le C_mh^{2},\eqno{(4.11)}$$
 $$\left\|\langle x\rangle^{-s}M_m^{(2)}(h,\theta)\langle x\rangle^{s}\right\|_{L^2\to L^2}\le C_mh^{2m+2},\eqno{(4.12)}$$
 which in turn follow from (3.19) and (3.20) (which clearly holds with $G$ replaced by $G_\theta$).
 \eproof
 
 Define the operators $\widetilde T_{j,\theta}^\flat$ by replacing in the definition of $\widetilde T_{j,\theta}$ the operator $Q_\theta(h)$ by $P_m^{(1)}(h,\theta)$. In precisely the same way as in the proof of (3.8) above, using (4.6) instead of (3.14), we get
 $$\left\|\widetilde T_{j,\theta}(t,h)-\widetilde T_{j,\theta}^\flat(t,h)\right\|_{L^1\to L^\infty}\le Ch^{2m+2+j-n/2} t^{-n/2}\le \widetilde Cht^{-n/2},\eqno{(4.13)}$$
provided $m$ is taken big enough. Therefore, it suffices to prove (4.2) with $\widetilde T_{j,\theta}$ replaced by $\widetilde T_{j,\theta}^\flat$. We will first do so for $j=1$. Let $0<\gamma\ll 1$ be a parameter to be fixed later on, depending on $h$. For $t\ge 2\gamma$, we have
 $$\left\|\widetilde T_{1,\theta}^\flat(t,h)\right\|_{L^1\to L^\infty}\le Ch^2\left\|\int_0^t  e^{i(t-\tau)G_0}\psi_1(h^2G_0)V_\theta P_m^{(1)}(h,\theta)e^{i\tau G_0}\psi_1(h^2G_0)d\tau\right\|_{L^1\to L^\infty}$$
 $$+C\left\|\int_0^t  e^{i(t-\tau)G_0}\psi_1(h^2G_0)V_\theta \left(P_m^{(1)}(h,\theta)-1\right)e^{i\tau G_0}\psi_1(h^2G_0)d\tau\right\|_{L^1\to L^\infty}$$
  $$\le Ct^{-n/2}\left(\int_0^\gamma+\int_{t-\gamma}^t\right)\left(\left\|e^{-i\tau G_0}V_\theta P_m^{(1)}(h,\theta) e^{i\tau G_0}\right\|_{L^1\to L^1}+\left\|e^{-i\tau G_0}V_\theta e^{i\tau G_0}\right\|_{L^1\to L^1}\right)d\tau$$
  $$+Ch^2\int_\gamma^{t/2}\left\|  e^{i(t-\tau)G_0}\psi_1(h^2G_0)\langle x\rangle^{-n/2-\epsilon'}\right\|_{L^2\to L^\infty}\left\|\langle x\rangle^{-n/2+1-\epsilon'}e^{i\tau G_0}\psi_1(h^2G_0)\right\|_{L^1\to L^2}d\tau$$
  $$+Ch^2\int_{t/2}^{t-\gamma}\left\|  e^{i(t-\tau)G_0}\psi_1(h^2G_0)\langle x\rangle^{-n/2+1-\epsilon'}\right\|_{L^2\to L^\infty}\left\|\langle x\rangle^{-n/2-\epsilon'}e^{i\tau G_0}\psi_1(h^2G_0)\right\|_{L^1\to L^2}d\tau$$
  $$\le C\gamma\theta^{k-n/2-\epsilon}t^{-n/2}+Cht^{-n/2}\int_\gamma^\infty\tau^{-n/2+1-\epsilon'}
  d\tau$$ $$= C\gamma\theta^{k-n/2-\epsilon}t^{-n/2}+C\gamma^{-n/2+2-\epsilon'}ht^{-n/2},\eqno{(4.14)}$$
  where we have used Lemma 4.2 together with (4.11) and (3.11) (with $j=0$). Clearly, (4.14) still holds for $0<t\le 2\gamma$. Choosing $\gamma$ such that $\gamma^{-n/2+2-\epsilon'}=h^{-1/2}$, we deduce the desired estimate from (4.14).
  
  Let now $j\ge 2$. Then $\widetilde T_{j,\theta}^\flat=F_{j,\theta}^\flat$, where the operators $F_{j,\theta}^\flat$, $j=0,1,...,$ are defined as follows
  $$F_{0,\theta}^\flat(t,h)=F_0(t,h)=e^{itG_0}\psi_1(h^2G_0),$$
$$F_{j,\theta}^\flat(t,h)=i\int_0^tF_0(t-\tau,h)V_\theta P_m^{(1)}(h,\theta)F_{j-1,\theta}^\flat(\tau,h)d\tau,\quad j\ge 1.$$
Let $0<\gamma\ll 1$ be a parameter to be fixed later on, depending on $h$. By Lemma 4.2, for 
$0<t\le 2j\gamma$, $j\ge 1$, we get
 $$\left\|F_{j,\theta}^\flat(t,h)\right\|_{L^1\to L^\infty}\le C\theta^{k-n/2-\epsilon}t^{-n/2}
 \int_0^{2j\gamma}\left\|e^{-i\tau G_0}F_{j-1,\theta}^\flat(\tau,h)\right\|_{L^1\to L^1}d\tau.\eqno{(4.15)}$$
 By induction in $j$, it is easy to see that we have the bound
  $$\left\|e^{-it G_0}F_{j,\theta}^\flat(t,h)\right\|_{L^1\to L^1}\le C_jt^j\theta^{j(k-n/2-\epsilon)},\quad\forall j\ge 0.\eqno{(4.16)}$$
  Clearly, (4.16) is trivial for $j=0$. Suppose that it holds for $j-1$. By Lemma 4.2 we have
 $$\left\|e^{-it G_0}F_{j,\theta}^\flat(t,h)\right\|_{L^1\to L^1}\le C\theta^{k-n/2-\epsilon}\int_0^t\left\|e^{-i\tau G_0}F_{j-1,\theta}^\flat(\tau,h)\right\|_{L^1\to L^1}d\tau$$
  $$\le C'_j\theta^{j(k-n/2-\epsilon)}\int_0^t\tau^{j-1}d\tau= C_jt^j\theta^{j(k-n/2-\epsilon)},$$
  which proves (4.16) for $j$. By (4.15) and (4.16), we conclude
 $$\left\|F_{j,\theta}^\flat(t,h)\right\|_{L^1\to L^\infty}\le C_j\gamma^j\theta^{j(k-n/2-\epsilon)}t^{-n/2},\quad 0<t\le 2j\gamma, \eqno{(4.17)}$$
  with a constant $C_j>0$ independent of $t$, $h$, $\theta$ and $\gamma$. We would like to obtain a similar estimate when $t\ge 2j\gamma$. To this end, decompose $F_{j,\theta}^\flat$ as follows
  $$F_{j,\theta}^\flat(t,h)=i\left(\int_0^\gamma +\int_\gamma^{t-\gamma}+\int_{t-\gamma}^t\right)F_0(t-\tau,h)V_\theta P_m^{(1)}(h,\theta)F_{j-1,\theta}^\flat(\tau,h)d\tau=:\sum_{\ell =1}^3 E_{j,\theta}^{(\ell)}(t,h,\gamma).$$
  By Lemma 4.2 and (4.16),
  $$\left\|E_{j,\theta}^{(1)}(t,h,\gamma)\right\|_{L^1\to L^\infty}\le C_j\gamma^j\theta^{j(k-n/2-\epsilon)}t^{-n/2}. \eqno{(4.18)}$$
 Clearly, the estimate (3.11) holds with $F_j$ replaced by $F_{j,\theta}^\flat$. Using this we obtain
 $$\left\|E_{j,\theta}^{(2)}(t,h,\gamma)\right\|_{L^1\to L^\infty}$$ $$\le
 C\int_\gamma^{t/2}\left\|F_0(t-\tau,h)\langle x\rangle^{-n/2-\epsilon'}\right\|_{L^2\to L^\infty}\left\|\langle x\rangle^{-n/2+1-\epsilon'}F_{j-1,\theta}^\flat(\tau,h)\right\|_{L^1\to L^2}d\tau$$
  $$+C\int_{t/2}^{t-\gamma}\left\|F_0(t-\tau,h)\langle x\rangle^{-n/2+1-\epsilon'}\right\|_{L^2\to L^\infty}\left\|\langle x\rangle^{-n/2-\epsilon'}F_{j-1,\theta}^\flat(\tau,h)\right\|_{L^1\to L^2}d\tau$$
  $$\le Ch^{j-2+\epsilon'/2}t^{-n/2}\int_\gamma^\infty\tau^{-n/2+1-\epsilon'}
  d\tau= C_j\gamma^{-n/2+2-\epsilon'}h^{j-2+\epsilon'/2}t^{-n/2},\eqno{(4.19)}$$
  with constants $C_j,\epsilon'>0$ independent of $t$, $h$, $\theta$ and $\gamma$ ($\epsilon'$ depending only on $\delta$). 
  Similarly, using (3.22), we also get
  $$\left\|\langle x\rangle^{-n/2-\epsilon'}E_{j,\theta}^{(2)}(t,h,\gamma)\right\|_{L^1\to L^2}\le C_j\gamma^{-n/2+2-\epsilon'}h^{n/2+j-2+\epsilon'/2}t^{-n/2}.\eqno{(4.20)}$$
 In what follows we will show that the operator $E_{j,\theta}^{(3)}$ satisfies the estimate
  $$\left\|E_{j,\theta}^{(3)}(t,h,\gamma)\right\|_{L^1\to L^\infty}\le C_j\gamma^{-n/2+2-\epsilon'}h^{j-2+\epsilon'/2}t^{-n/2}+C_j\gamma^j\theta^{j(k-n/2-\epsilon)}t^{-n/2},\quad  t\ge 2j\gamma. \eqno{(4.21)}$$ 
  To this end, it suffices to show that modulo operators satisfying (4.19), the operator $E_{j,\theta}^{(3)}$ is a finite sum of operators of the form
  $$\int_{I_1}\int_{I_2}...\int_{I_j}F_0(t-\tau_1,h)V_\theta P_m^{(1)}(h,\theta)F_0(\tau_1-\tau_2,h)V_\theta P_m^{(1)}(h,\theta)...V_\theta P_m^{(1)}(h,\theta)F_0(\tau_j,h)d\tau_1d\tau_2...d\tau_j,\eqno{(4.22)}$$
 where $I_\nu$, $\nu=1,...,j,$ are intervals of length $|I_\nu|=O(\gamma)$. Indeed, by Lemma 4.2 an operator of form (4.22) satisfies (4.18).
 We will show that given any integer $1\le\nu\le j-1$, the operator $E_{j,\theta}^{(3)}$ can be written in the form
 $$\int_0^\gamma...\int_0^\gamma F_0(\tau'_1,h)V_\theta P_m^{(1)}(h,\theta)...F_0(\tau'_\nu,h)V_\theta P_m^{(1)}(h,\theta)E_{j-\nu,\theta}^{(3)}(t-\tau'_1-...-\tau'_\nu,h,\gamma)d\tau'_1...d\tau'_\nu,\eqno{(4.23)}$$ 
 modulo operators satisfying (4.19) and operators of the form (4.22). This would imply the desired result because the operator (4.23) with $\nu=j-1$ is of the form (4.22). We will proceed by induction in $\nu$. Let us see that the claim holds true for $\nu=1$. We write 
 $$E_{j,\theta}^{(3)}(t,h,\gamma)=\sum_{\ell =1}^3\int_0^\gamma F_0(\tau'_1,h)V_\theta P_m^{(1)}(h,\theta)E^{(\ell)}_{j-1,\theta}(t-\tau'_1,h)d\tau'_1.\eqno{(4.24)}$$ 
 Clearly, the first operator in the sum in the right-hand side of (4.24) is of the form (4.22), while the third one is of the form (4.23) with $\nu=1$. On the other hand, using (3.11) with $j=0$, $s=1/2-\epsilon/2$, together with (4.20), it is easy to see that the second one satisfies (4.19). Suppose now that the claim holds true for some $\nu$, $1\le\nu\le j-2$. Then we decompose the operator in (4.23) as follows
 $$\sum_{\ell =1}^3\int_0^\gamma...\int_0^\gamma F_0(\tau'_1,h)V_\theta P_m^{(1)}(h,\theta)...F_0(\tau'_{\nu+1},h)V_\theta P_m^{(1)}(h,\theta)$$ $$\times E_{j-\nu-1,\theta}^{(\ell)}(t-\tau'_1-...-\tau'_{\nu+1},h,\gamma)d\tau'_1...d\tau'_{\nu+1}.\eqno{(4.25)}$$ 
 Clearly, the first operator in the sum in (4.25) is of the form (4.22), while the third one is of the form (4.23) with $\nu+1$. Therefore, to prove the claim it sufices to show that the second one satisfies (4.19). However, this follows easily from (4.20) and the following consequence of (3.11) with $j=0$, $s=1/2-\epsilon/2$, and (3.22):
 $$\left\|F_0(\tau'_1,h)V_\theta P_m^{(1)}(h,\theta)...F_0(\tau'_{\nu+1},h)\langle x\rangle^{-1-\epsilon}\right\|_{L^2\to L^\infty}$$ $$\le C_{\epsilon,\nu}h^{1-n/2-\epsilon/2}(\tau'_1)^{-1+\epsilon/2}\langle\tau'_2/h\rangle^{-1-\epsilon}...\langle\tau'_{\nu+1}/h\rangle^{-1-\epsilon},\eqno{(4.26)}$$ 
 for every $0<\epsilon\ll 1$.
  
  By (4.17), (4.18), (4.19) and (4.21), we conclude that the operators $F_{j,\theta}^\flat$, $j\ge 2$, satisfy the estimate
  $$\left\|F_{j,\theta}^\flat(t,h)\right\|_{L^1\to L^\infty}\le C_j\gamma^{-n/2+2-\epsilon'}h^{j-2+\epsilon'/2}t^{-n/2}+C_j\gamma^j\theta^{j(k-n/2-\epsilon)}t^{-n/2},\quad \forall t>0. \eqno{(4.27)}$$
  Choosing $\gamma$ such that
  $$\gamma^{-n/2+2-\epsilon'}=h^{-j+2-\epsilon'/4},$$
 we deduce the desired estimate from (4.27).
 \eproof
 
 \section{Proof of Theorem 1.2}
 
 It is easy to see that Theorem 1.2 follows from the following
 
 \begin{Theorem} Under the assumptions of Theorem 1.2, there exist an integer $m$ and constants $C,\beta>0$ such that the following dispersive estimates hold true for all $0<h\ll 1$, $t>0$:
$$\|{\cal F}_j(t)\|_{L^1\to L^\infty}\le Ct^{-n/2},\quad 1\le j\le m,\eqno{(5.1)}$$
$$\left\|e^{itG}\psi(h^2G)-\sum_{j=0}^m{\cal F}_j(t)\psi(h^2G)\right\|_{L^1\to L^\infty}\le Ch^\beta t^{-n/2}.\eqno{(5.2)}$$
\end{Theorem}

{\it Proof.} Let us first see that (5.1) holds for all $j\ge 0$. It is trivial for $j=0$, while when $0<t\le 2$ it follows from (1.4). Let now $t\ge 2$. We will proceed by induction in $j$. Suppose that (5.1) holds for $j-1$. This implies
$$\left\|\int_1^{t-1}{\cal F}_{j-1}(t-\tau)V{\cal F}_0(\tau)d\tau\right\|_{L^1\to L^\infty}\le \|V\|_{L^1}\int_1^{t-1}\|{\cal F}_{j-1}(t-\tau)\|_{L^1\to L^\infty}\|{\cal F}_0(\tau)\|_{L^1\to L^\infty}d\tau$$
$$\le C\int_1^{t-1}(t-\tau)^{-n/2}\tau^{-n/2}d\tau\le C't^{-n/2}.\eqno{(5.3)}$$
Clearly, (5.1) for $j$ follows from (5.3) and (1.5) applied with $\gamma=1$, $m_1=j-1$, $m_2=0$.

Let $0<\gamma\le 1$ be a parameter to be fixed later on, depending on $h$. Set $\gamma_1=\gamma$ if $t\ge 2\gamma$, $\gamma_1=t/2$ if $t\le 2\gamma$. Iterating Duhamel's formula we obtain the identity
$$e^{itG}=\sum_{j=0}^{m_1}{\cal F}_j(t)+i\int_0^t{\cal F}_{m_1}(t-\tau)V\left(e^{i\tau G}-e^{i\tau G_0}\right)d\tau$$
$$=\sum_{j=0}^{m_1}{\cal F}_j(t)+i\int_{\gamma_1}^t{\cal F}_{m_1}(t-\tau)V\left(e^{i\tau G}-e^{i\tau G_0}\right)d\tau$$
$$+i\sum_{\nu=1}^{m_2-1}\int_0^{\gamma_1}{\cal F}_{m_1}(t-\tau)V{\cal F}_\nu(\tau)d\tau+i^2\int_0^{\gamma_1}\int_0^\tau 
{\cal F}_{m_1}(t-\tau)V{\cal F}_{m_2}(\tau -s)Ve^{isG}ds d\tau,\eqno{(5.4)}$$
for all integers $m_1,m_2\ge 2$. Hence
$$\left\|e^{itG}\psi(h^2G)-\sum_{j=0}^m{\cal F}_j(t)\psi(h^2G)\right\|_{L^1\to L^\infty}$$ $$\le \int_{\gamma_1}^t\left\|{\cal F}_{m_1}(t-\tau)\right\|_{L^1\to L^\infty}\left\|V\left(e^{i\tau G}-e^{i\tau G_0}\right)\psi(h^2G)\right\|_{L^1\to L^1}d\tau$$
 $$+C\sum_{\nu=1}^{m_2-1}\left\|\int_0^{\gamma_1}{\cal F}_{m_1}(t-\tau)V{\cal F}_\nu(\tau)d\tau\right\|_{L^1\to L^\infty}$$
  $$+\int_0^{\gamma_1}\int_0^\tau 
\left\|{\cal F}_{m_1}(t-\tau)\right\|_{L^1\to L^\infty}\left\|V{\cal F}_{m_2}(\tau -s)\right\|_{L^1\to L^1}\left\|Ve^{isG}\psi(h^2G)\right\|_{L^1\to L^1}ds d\tau$$
 $$\le C\int_{\gamma_1}^t\left\|{\cal F}_{m_1}(t-\tau)\right\|_{L^1\to L^\infty}\left\|\langle x\rangle^{-n/2-\epsilon'}\left(e^{i\tau G}-e^{i\tau G_0}\right)\psi(h^2G)\right\|_{L^1\to L^2}d\tau$$
 $$+C\sum_{\nu=1}^{m_2-1}\left\|\int_0^{\gamma_1}{\cal F}_{m_1}(t-\tau)V{\cal F}_\nu(\tau)d\tau\right\|_{L^1\to L^\infty}$$
  $$+C\int_0^{\gamma_1}\int_0^\tau 
\left\|{\cal F}_{m_1}(t-\tau)\right\|_{L^1\to L^\infty}\left\|{\cal F}_{m_2}(\tau -s)\right\|_{L^1\to L^\infty}\left\|\langle x\rangle^{-n/2-\epsilon'}e^{isG}\psi(h^2G)\right\|_{L^1\to L^2}ds d\tau.\eqno{(5.5)}$$
On the other hand, it is 
proved in \cite{kn:V2} (see Proposition 4.1) for potentials satisfying (1.2) with $\delta>\frac{n+2}{2}$ that we have the estimate
$$\left\|\langle x\rangle^{-n/2-\epsilon'}\left(e^{itG}\psi(h^2G)-e^{itG_0}\psi(h^2G_0)\right)\right\|_{L^1\to L^2}\le Ch^{1-\epsilon}t^{-n/2}.\eqno{(5.6)}$$
Hence
$$\left\|\langle x\rangle^{-n/2-\epsilon'}\left(e^{itG}-e^{itG_0}\right)\psi(h^2G)\right\|_{L^1\to L^2}$$ $$\le \left\|\langle x\rangle^{-n/2-\epsilon'}\left(e^{itG}\psi(h^2G)-e^{itG_0}\psi(h^2G_0)\right)\right\|_{L^1\to L^2}$$ $$+C\left\|\psi(h^2G)-\psi(h^2G_0)\right\|_{L^1\to L^1}\left\|e^{itG_0}\right\|_{L^1\to L^\infty}$$
 $$\le Ch^{1-\epsilon}t^{-n/2}+Ch^2t^{-n/2}\le C'h^{1-\epsilon}t^{-n/2}.\eqno{(5.7)}$$
 Using (5.7) together with (1.4) and (5.1) we can bound the first integral in the right-hand side of (5.5) by
 $$Ch^{1-\epsilon}\left(\int_{\gamma_1}^{t/2}+\int_{t/2}^t\right)\left\|{\cal F}_{m_1}(t-\tau)\right\|_{L^1\to L^\infty}\tau^{-n/2}d\tau$$
  $$\le Ch^{1-\epsilon}t^{-n/2}\int_{\gamma_1}^{t/2}\tau^{-n/2}d\tau+Ch^{1-\epsilon}t^{-n/2}\int_{t/2}^t\left\|{\cal F}_{m_1}(t-\tau)\right\|_{L^1\to L^\infty} d\tau$$
   $$\le Ch^{1-\epsilon}t^{-n/2}\int_{\gamma}^\infty\tau^{-n/2}d\tau+Ch^{1-\epsilon}t^{-n/2}\int_0^\infty\left\|{\cal F}_{m_1}(\tau')\right\|_{L^1\to L^\infty} d\tau'$$
 $$\le Ch^{1-\epsilon}\gamma^{-(n-2)/2}t^{-n/2}+Ch^{1-\epsilon}t^{-n/2}\le  C'h^{1-\epsilon}\gamma^{-(n-2)/2}t^{-n/2},$$
 provided $m_1$ is taken big enough. Furthermore, in view of (1.5), each term in the sum in the right-hand side of (5.5) is bounded by
 $C\gamma^\varepsilon t^{-n/2}$ for all $t>0$.  To bound the last integral we will use the following estimate proved in \cite{kn:V2} (see Propositions 2.1 and 4.1)
$$\left\|\langle x\rangle^{-s-1/2-\epsilon}e^{it G}\psi(h^2G)\right\|_{L^1\to L^2}\le Ch^{s-(n-1)/2}|t|^{-s-1/2}\eqno{(5.8)}$$
for every $0<\epsilon\ll 1$, $1/2-\epsilon/4\le s\le (n-1)/2$, $0<h\ll 1$, $t\neq 0$. 
Using (5.8) together with (1.4) and (5.1) we bound the integral under question by
$$Ch^{-(n-2)/2-\epsilon/4}t^{-n/2}\int_0^{\gamma}\int_0^\tau (\tau-s)^{\varepsilon m_2-n/2}s^{-1+\epsilon/4}dsd\tau$$
 $$\le Ch^{-(n-2)/2-\epsilon/4}t^{-n/2}\int_0^{\gamma}\tau^{\varepsilon m_2-n/2}d\tau\le C\gamma^{\varepsilon m_2-n/2+1}h^{-(n-2)/2-\epsilon/4}t^{-n/2},$$
provided $m_2$ is taken big enough. Summing up the above estimates, we conclude that the left-hand side of (5.5) is bounded by
$$Ch^{1-\epsilon}\gamma^{-(n-2)/2}t^{-n/2}+C\gamma^\varepsilon t^{-n/2}+C\gamma^{\varepsilon m_2-n/2+1}h^{-(n-2)/2-\epsilon/4}t^{-n/2}\eqno{(5.9)}$$
for all $t>0$ and all $0<\gamma\le 1$. Take $\gamma=h^{1/(n-2)}$ and fix $m_2$ so that
$$\gamma^{\varepsilon m_2-n/2+1}h^{-(n-2)/2}\le h.$$
Hence (5.9) is $O(h^\beta)t^{-n/2}$ with some $\beta>0$, which is the desired result.
\eproof

 {\bf Acknowledgements.} We would like to thank William Green for bringing to our attention the papers \cite{kn:EG}, \cite{kn:FY}, \cite{kn:Y}. The first two authors have been partially supported by the CNPq-Brazil.

F. Cardoso

Universidade Federal de Pernambuco, 

Departamento de Matem\'atica, 

CEP. 50540-740 Recife-Pe, Brazil,

e-mail: fernando@dmat.ufpe.br\\

C. Cuevas

Universidade Federal de Pernambuco, 

Departamento de Matem\'atica, 

CEP. 50540-740 Recife-Pe, Brazil,

e-mail: cch@dmat.ufpe.br\\

G. Vodev

Universit\'e de Nantes,

 D\'epartement de Math\'ematiques, UMR 6629 du CNRS,
 
 2, rue de la Houssini\`ere, BP 92208, 
 
 44332 Nantes Cedex 03, France,
 
 e-mail: vodev@math.univ-nantes.fr

\end{document}